\documentclass[12pt]{amsart}

\usepackage{amscd}
\usepackage{mathabx}
\input bookman.sty
\boldmath

\usepackage{comment}
\excludecomment{excl}

\usepackage{xcolor}

\usepackage{helvet}

\usepackage{graphics}
\usepackage{amssymb}
\usepackage{amsxtra}
\usepackage{amsmath}
\usepackage{mathrsfs}

\usepackage[arrow, matrix]{xy}
\xyoption{frame}

\theoremstyle{plain}

\newtheorem{theo}{Theorem}[section]
\newtheorem{lemm}[theo]{Lemma}
\newtheorem{prop}[theo]{Proposition}

\newtheorem{coro}[theo]{Corollary}

\theoremstyle{definition}

\newtheorem{defi}[theo]{Definition}
\newtheorem{rema}[theo]{Remark}

\newfont{\rmm}{cmr10 scaled 1000}

\newfont{\itt}{cmsl10 scaled 1000}

\newfont{\rM}{cmr10 scaled 1700}

\renewcommand{\frak}{\mathfrak}

\newcounter{lemma}[section]

\newcounter{tempcounter}

\newcommand{\lb}{\label}

\newcommand{\rrf}[1]{(\ref{#1})}

\renewcommand{\a}{\alpha}
\renewcommand{\b}{\beta}

\newcommand{\ve}{\varepsilon}

\renewcommand{\t}{\theta}
\renewcommand{\l}{\lambda}

\renewcommand{\r}{\rho}

\newcommand{\s}{\sigma}

\renewcommand{\o}{\omega}

\renewcommand{\O}{\Omega}

\newcommand{\II}{{\mathcal I}}

\newcommand{\cc}{{\mathbb{C}}}

\newcommand{\hh}{{\mathbb{H}}}

\newcommand{\nn}{{\mathbb{N}}}

\newcommand{\qq}{{\mathbb{Q}}}
\newcommand{\rr}{{\mathbb{R}}}

\newcommand{\ttt}{{\mathbb{T}}}

\newcommand{\zz}{{\mathbb{Z}}}

\newcommand{\OOOO}{{\mathscr{O}}}
\newcommand{\PPPP}{{\mathscr{P}}}

\newcommand{\TTTT}{{\mathscr{T}}}

\renewcommand{\Im}{\text{\rm Im }}

\newcommand{\Id}{\text{\rm Id}}

\newcommand{\bere}{\begin{rema}}
\newcommand{\bede}{\begin{defi}}

\renewcommand{\beth}{\begin{theo}}
\newcommand{\bele}{\begin{lemm}}
\newcommand{\bepr}{\begin{prop}}
\newcommand{\beeq}{\begin{equation}}
\newcommand{\bega}{\begin{gather}}
\newcommand{\begaa}{\begin{gather*}}
\newcommand{\been}{\begin{enumerate}}

\newcommand{\bedee}{\begin{defii}}
\newcommand{\bethh}{\begin{theoo}}
\newcommand{\belee}{\begin{lemmm}}
\newcommand{\beprr}{\begin{propp}}

\newcommand{\beco}{\begin{coro}}

\newcommand{\beal}{\begin{aligned}}

\newcommand{\enre}{\end{rema}}

\newcommand{\enco}{\end{coro}}
\newcommand{\enpr}{\end{prop}}
\newcommand{\enth}{\end{theo}}
\newcommand{\enle}{\end{lemm}}
\newcommand{\enen}{\end{enumerate}}
\newcommand{\enga}{\end{gather}}
\newcommand{\engaa}{\end{gather*}}
\newcommand{\eneq}{\end{equation}}
\newcommand{\enal}{\end{aligned}}

\newcommand{\bq}{\begin{equation}}
\newcommand{\bqq}{\begin{equation*}}

\renewcommand{\leq}{\leqslant}
\renewcommand{\geq}{\geqslant}

\newcommand{\ove}{\overline}

\newcommand{\sbs}{\subset}

\newcommand{\sut}{~such~that~}

\newcommand{\wrt}{with respect to}
\newcommand{\ho}{homomorphism}

\newcommand{\ma}{manifold}

\newcommand{\Prf}{{\it Proof.\quad}}

\newcommand{\chart}{\Phi_p:U_p\to B^n(0,r_p)}
\newcommand{\atlas}{\{\Phi_p:U_p\to B^n(0,r_p)\}_{p\in S(f)}}

\newcommand{\qs}{\hfill\square}

\newcommand{\pa}{\vskip0.1in}

\newcommand{\arrh}[3]
{
\xymatrix{
{#1} \ar[r]^<<<<{#2}  &{#3}
}
}

\newcommand{\arrr}[1]
{\arrh {}{#1}{}}

\newcommand{\arr}
{\arrr {}}

\newcommand{\arrto}
{\xymatrix{{} \ar@{|-{>}}[r]  & {} } }

\newcommand{\arrinto}
{\xymatrix{{} \ar@{^{(}->}[r]  & {} } }

\newcommand{\OT}{K. Oeljeklaus\ and\ M. Toma}

\begin{document}

\title
[On generalized Inoue manifolds]
{On generalized Inoue manifolds}
\author{Hisaaki Endo  and  Andrei Pajitnov}
\address{Department of Mathematics
Tokyo Institute of Technology
2-12-1 Ookayama, Meguro-ku
Tokyo 152-8551
Japan}
\email{endo@math.titech.ac.jp}
\address{Laboratoire Math\'ematiques Jean Leray 
UMR 6629,
Universit\'e de Nantes,
Facult\'e des Sciences,
2, rue de la Houssini\`ere,
44072, Nantes, Cedex}                    
\email{andrei.pajitnov@univ-nantes.fr}

\thanks{} 
\begin{abstract}
This paper is about a generalization of famous Inoue's surfaces.
Let $M$ be a matrix in $SL(2n+1,\mathbb{Z})$
having only one real eigenvalue which is simple.
We associate to $M$ 
a complex manifold $T_M$ of complex dimension $n+1$. 
This manifold 
fibers over $S^1$ with the fiber $\mathbb{T}^{2n+1}$ and monodromy
$M^\top$. Our construction is elementary and does not use 
algebraic number theory. We show that some of the Oeljeklaus-Toma manifolds
are biholomorphic to the manifolds of type $T_M$.
We prove that if $M$ is not diagonalizable, then $T_M$ 
does not admit a K\"ahler structure and is not
homeomorphic to any of Oeljeklaus-Toma manifolds.
\end{abstract}
\maketitle
\tableofcontents
\newcommand{\ssll}{\mathrm{SL}}

\section{Introduction}
\label{s:intro}
\subsection{Background}
\label{su:bckgr}

In 1972 M. Inoue \cite{Inoue1974}
constructed complex surfaces
having remarkable properties:
they have second Betti number equal to zero
and contain no complex curves.
These surfaces (called  now {\it Inoue surfaces})
attracted a lot of attention.
It was proved by F. Bogomolov \cite{Bogomolov}
(see also A. Teleman \cite{T}) 
that each complex surface
with $b_2(X)=0$ containing no complex
curves is isomorphic to an Inoue surface.
Inoue surfaces are not algebraic, and moreover 
they do not admit K\"ahler metric
(since their first Betti number is odd).

Let us say that a matrix $M\in \ssll(2n+1,\zz)$ 
is of  type $\II$, if it has only one real 
eigenvalue which is irrational and simple.
Inoue's construction associates to every such 
matrix $M\in \ssll(3,\zz)$
a complex surface $T_M$ obtained as a quotient 
of $\hh\times \cc$ by action of a discrete group
(here $\hh$ is the upper half-plane).
This manifold fibers over  $S^1$ with fiber $\ttt^3$
and the monodromy of this fibration equals 
the diffeomorphism of $\ttt^{3}$ determined by 
$M^{\top}$.

Inoue's construction was generalized
to higher dimensions in several papers
in particular in 
a celebrated paper of  K. Oeljeklaus and M. Toma 
\cite{OT}.
The construction of Oeljeklaus and Toma
uses algebraic number theory.
It starts with an algebraic number field $K$.
Denote by $s$ the number of  embeddings of $K$ to $\rr$
and by $2t$  number of  non-real embeddings of $K$ to $\cc$,
so that $(K:\qq)=s+2t$.
\OT~ constructed an action of a certain semi-direct 
product $\zz^s \ltimes \zz^{2t+s}$
on $\hh^s\times \cc^t$, \sut~
the quotient is a compact complex manifold
of complex dimension $s+t$.
This \ma~  has interesting geometric properties, studied in details in \cite{OT};
in particular, it does not admit a K\"ahler metric.
The original Inoue surface corresponds to the 
algebraic number field generated by the eigenvalues
of the matrix $M$.

\subsection{Outline of the paper}
\label{su:outline}

In the present paper we introduce another generalization
of Inoue's construction.
Our method does not use algebraic number theory,
it generalizes the original Inoue's approach.
Let  and
$M\in \ssll(2n+1,\zz)$ 
be a matrix of  type $\II$.
We construct an action of a certain 
semi-direct product $\zz\ltimes \zz^{2n+1}$
on $\hh\times \cc^n$, the quotient is a complex
non-K\"ahler manifold $T_M$.
It fibers over  $S^1$ with fiber $\ttt^{2n+1}$
and the monodromy of this fibration equals 
the diffeomorphism of $\ttt^{2n+1}$ determined by 
$M^{\top}$.
The construction of the manifolds is done in 
Section \ref{s:constr}
and their properties are studied
in Sections 
\ref{s:topo_pr}
 and 
 \ref{s:geom_pr}.

The basic difference of our construction
from the preceding generalizations of Inoue's
work is that the matrix $M$ 
can be non-diagonalizable. In the  non-diagonalizable case 
the monodromy of the fibration $T_M\to S^1$ is also non-diagonalizable,
and this implies that  certain Massey products
in the cohomology of $T_M$ with local coefficients
do not vanish. Therefore  the manifold
$T_M$ is not strongly formal
(see \cite{KP} for definition of strongly formal manifolds).
This implies that  the manifold $T_M$ and its cartesian powers $(T_M)^l, \ l\in\nn$
do not admit a structure of a K\"ahler manifold.

In Section \ref{s:Tm-OT}
we show that some of the Oeljeklaus-Toma
manifolds are biholomorphic to manifolds $T_M$ 
for some special choices of the matrix $M$.
Then  we show that if 
$M$ is non-diagonalizable  then the manifold 
$T_M$ is not homeomorphic to any of 
Oeljeklaus-Toma
manifolds 
(see Subsection 
\ref{su:non-diag}).

\section{Manifold $T_M$: the construction}
\lb{s:constr}

\bigskip

Let $n$ be a positive integer and 
$M=(m_{ij})$ an element of $\ssll(2n+1,\mathbb{Z})$. 
Suppose that $M$ has eigenvalues 
$\alpha,\beta_1,\ldots ,\beta_k,\overline{\beta}_1, \ldots ,\overline{\beta}_k \in\mathbb{C}$ 
which satisfy $\alpha\in\mathbb{R}, \alpha>0, \alpha\ne 1, \frak{Im}(\beta_j)>0$ 
for every $j\in\{1,\ldots ,k\}$, 
and the eigenspace $V$ of $M$ corresponding to $\alpha$ has dimension one. 
We denote the generalized eigenspace of $M$ corresponding 
to an eigenvalue $\beta$ 
by $W(\beta)$, namely 
\[
W(\beta):=\{x\in\mathbb{C}^{2n+1}\, |\, (M-\beta I)^Nx=0 \; 
{\rm for \ some \ positive \ integer}\; N\}.
\]
We obtain a direct sum decomposition of $\mathbb{C}^{2n+1}$ 
into complex $M$-invariant subspaces 
\[
\mathbb{C}^{2n+1}
=
V\oplus W\oplus \overline{W}, 
\quad 
W:=\bigoplus_{j=1}^k W(\beta_j), 
\quad 
\overline{W}:=\bigoplus_{j=1}^k W(\overline{\beta}_j).
\]

Let $a$ be a real eigenvector of $M$ corresponding to $\alpha$ and 
$b_1,\ldots ,b_n$ a basis of $W$. 
Then $\overline{b}_1,\ldots ,\overline{b}_n$ is a basis of $\overline{W}$, 
and $a, b_1,\ldots ,b_n, \overline{b}_1,\ldots ,\overline{b}_n$ is a basis of $\mathbb{C}^{2n+1}$. 
Let $f_M:W\rightarrow W$ be the restriction of $M$ to $W$, 
and $R=(r_{ij})$ the matrix of $f_M$ in
the basis $b_1,\ldots ,b_n$, 
namely
\bq\lb{f:M}
\qquad 
Mb_j=\sum_{\ell=1}^n r_{\ell j}b_{\ell},
\quad 
(
r_{\ell j}\in\mathbb{C}). 
\end{equation}
Write
\[
a=
\left(
\begin{array}{c}
a^{(1)} \\
\vdots \\
a^{(2n+1)} 
\end{array}
\right), \ \ 
\; 
b_1=
\left(
\begin{array}{c}
b_1^{(1)} \\
\vdots \\
b_1^{(2n+1)} 
\end{array}
\right), 
\ldots ,
\; 
b_n=
\left(
\begin{array}{c}
b_n^{(1)} \\
\vdots \\
b_n^{(2n+1)} 
\end{array}
\right),
\]
where $a^{(1)},\ldots ,a^{(2n+1)}$ are real numbers. 
We also consider the vectors
$$
v_i:=
(a^{(i)},b_1^{(i)},\ldots ,b_n^{(i)})
\in\mathbb{R}\times\mathbb{C}_n=\mathbb{R}_{2n+1}
\ \ 
$$
$$
u_i:={v_i}^{\top} \in \rr\times \cc^n=\rr^{2n+1}, 
\quad 
(i\in\{1,\ldots ,2n+1\}).
$$

\bigskip
The following lemma is easy to prove.
\noindent
\bele\lb{l:lemmaA}
The vectors $v_1,\ldots ,v_{2n+1}$ are linearly independent over $\mathbb{R}$. 
$\qs$ \enle
\bigskip

We consider the following matrices and vectors: 
\begin{align*}
B
& :=
\left(
b_1, \ldots ,b_n
\right)
=
\left(
\begin{array}{ccc}
b_1^{(1)} & \cdots & b_n^{(1)} \\
\vdots & & \vdots \\
b_1^{(2n+1)} & \cdots & b_n^{(2n+1)}
\end{array}
\right), \\
b^{(i)}
& :=
\left(
b_1^{(i)}, 
\ldots , 
b_n^{(i)}
\right) \in \mathbb{C}_n 
\quad 
(i\in\{1,\ldots ,2n+1\}).
\end{align*}

\bigskip
A direct computation proves the following lemma.
\noindent
\bele\lb{l:LemmaB} \ 
The equality $MB=BR$ holds. 
In particular
$$
b^{(i)}R=\sum_{j=1}^{2n+1} m_{ij}b^{(j)}
\ \  
{\rm for\ every \ \ }
i\in\{1,\ldots ,2n+1\}. \quad\quad\qs$$
\enle
\bigskip

\noindent

Let $\mathbb{H}$ be the upper half of the complex plane, 
namely
\[
\mathbb{H}=\{w\in\mathbb{C}\, |\, \frak{Im}(w)>0\}.
\]
We consider complex-analytic automorphisms  
$g_0,g_1,\ldots ,g_{2n+1}:\mathbb{H}\times\mathbb{C}^n\rightarrow 
\mathbb{H}\times\mathbb{C}^n$ defined by
\[
g_0(w,z):=(\alpha w,\, R^{\top}z), 
\quad 
g_i(w,z):=
(w,z)+u_i
\]
for every $(w,z)\in \mathbb{H}\times\mathbb{C}^n$ and $i\in\{1,\ldots ,2n+1\}$. 
Let  $G_M$ be the subgroup of $\mathrm{Aut}(\mathbb{H}\times\mathbb{C}^n)$ 
generated by $g_0,g_1,\ldots ,g_{2n+1}$, 
$H_M$ the subgroup of $\mathrm{Aut}(\mathbb{H}\times\mathbb{C}^n)$ 
generated by $g_1,\ldots ,g_{2n+1}$, 
and $\langle g_0\rangle$ the infinite cyclic group generated by $g_0$. 
Lemma \ref{l:lemmaA}
implies that
$H_M$ is a free abelian group of rank $2n+1$. 

\bigskip

\noindent
\bele\lb{l:LemmaC} \ 
For every $i\geq 1$ we have
$$g_0g_ig_0^{-1}=g_1^{m_{i1}}\cdots g_{2n+1}^{m_{i,2n+1}}.$$
In particular, 
$H_M$ is a normal subgroup of $G_M$. 
\enle
{\it Proof.} 
Let $(w,z)\in\hh\times\cc$.
We have
$$
g_0(g_i(w,z))
=
g_0((w,z)+ u_i)
=
\left(\alpha(w+a^{(i)}), R^{\top}(z+(b^{(i)})^{\top})\right)
$$
$$
=
\left(\alpha w+\sum_{j=1}^{2n+1} m_{ij}a^{(j)}, \ \ 
R^{\top}z+\sum_{j=1}^{2n+1} m_{ij}(b^{(j)})^{\top} \right),\ \ 
$$
which \ equals \ \ 
$(g_1^{m_{i1}}\cdots g_{2n+1}^{m_{i,2n+1}})(g_0(w,z)) 
$
\ \ 
by Lemma \ref{l:LemmaB}. $\qs$

\bigskip
Observe that the group $G_M/H_M$ is generated by one element $g_0$.
For $(w,z)\in \hh\times \cc$ denote $\Im w$ by $p_1(w,z)$.
We have then 
$p_1(g_i(w,z))=p_1(w,z)$ for $i>0$
and 
$p_1(g_0(w,z))=\a\cdot p_1(w,z)$.
Therefore the element $g_0^n$ is not in $H_M$ for any $n\in\zz$.
\bepr\lb{p:PropE}
The group $G_M$ is isomorphic to a semi-direct product
of $\zz$ and $\zz^{2n+1}$
associated to the action of $\zz$ on $\zz^{2n+1}$
given by the formula
$t\cdot v= M^{\top} v$
(where $t$ is a generator of $\zz$ and $v\in\zz^{2n+1}$).
\enpr
\Prf
It follows from the observation above that the group $G_M/H_M$
is infinite cyclic. Therefore the exact sequence 
$$
1\longrightarrow\ H_M
\hookrightarrow
G_M \
{\longrightarrow} \ \
G_M/H_M \longrightarrow 1
$$
is isomorphic to
$$
1\longrightarrow \zz^{2n+1}
{\longrightarrow} \ \ G_M
{\longrightarrow} \ \
\zz \longrightarrow 1.
$$
The action of
the group $\zz$ on $\zz^{2n+1} $
is easily deduced from Lemma \ref{l:LemmaC}. $\qs$

\noindent
\beco\lb{c:CorollaryF} \ 
The group $G_M$ admits a finite presentation with generators 
$g_0,g_1,\ldots ,g_{2n+1}$ and defining relations 
\begin{align*} 
& g_ig_j=g_jg_i 
\quad 
(i,j\in\{1,\ldots ,2n+1\}), \\
& g_0g_ig_0^{-1}=g_1^{m_{i1}}\cdots g_{2n+1}^{m_{i,2n+1}} 
\quad 
(i\in\{1,\ldots ,2n+1\}). 
\end{align*}
\enco
\noindent
{\it Proof.} \ 
It follows from Lemma \ref{l:LemmaC} and 
Proposition \ref{p:PropE} 
(see \cite[Section 5.4]{Johnson1990}). 
$\qs$

\bigskip

\beco\lb{c:LemmaL}
The group $H_M$ includes the commutator subgroup $[G_M,G_M]$ of $G_M$, and 
the quotient $H_M/[G_M,G_M]$ is finite. 
\enco
\noindent
{\it Proof.} \ 
The first part of the Lemma follows
from the fact that $G_M/H_M$ is abelian.
Further, Corollary \ref{c:CorollaryF} implies that 
the group $H_M/[G_M,G_M]$ is isomorphic to the abelian group 
generated by $g_1,\ldots ,g_{2n+1}$ with relations 
\[
g_i=m_{i1}g_1+\cdots +m_{i,2n+1}g_{2n+1} 
\quad 
(i\in\{1,\ldots ,2n+1\}).
\]
Since $M$ does not have eigenvalue $1$, 
we see $\det(M-I)\ne 0$. 
Thus the group $H_M/[G_M,G_M]$ is finite. $\qs$

\noindent
\bepr\lb{p:PropG} \ 
The action of $G_M$ on $\mathbb{H}\times\mathbb{C}^n$ is free and 
properly discontinuous. 
\enpr
\newcommand{\HC}{\hh\times\cc^n}
\pa
\noindent
{\it Proof.} \ 
We will prove that the action is free,
the proof of the discontinuity is similar.
Let $(w,z)\in \hh\times\cc^n$, and $g\in G_M$.
Assume that 
$g(w,z)=(w,z)$ 
for some $(w,z)\in \HC$.
Write $g=g_0^{m_0}\cdot h$, where $h\in H_M$.
Observe that 
$p_1(g(w,z))=\alpha^{m_0}\cdot \frak{Im} w$; therefore
$m_0=0$, and $g\in H_M$.
THe action of $H_M$ leaves invariant the $(2n+1)$-dimensional
real affine subspace 
$$
V=\{(w',z')~|~ \frak{Im} w'=\frak{Im} w\}.$$
On this space $H_M$ acts
as a full lattice generated by vectors
$v_1, \ldots, v_{2n+1}$.
This action is free, therefore $g=1$. $\qs$

\bigskip

We consider the map 
$g_M:\mathbb{R}\times\mathbb{C}^n\rightarrow\mathbb{R}\times\mathbb{C}^n$ 
defined by 
\[
g_M(x,z):=(\alpha x, \, R^{\top}z) 
\quad 
(x\in\mathbb{R},\, z\in\mathbb{C}^n).
\]

\noindent
A direct computation using Lemma \ref{l:LemmaB}
proves the following Lemma.
\bele\lb{l:LemmaH} \ 
The matrix of the linear transformation $g_M$ with respect to the basis 
$(u_1, \ldots, u_{2n+1})$
is equal to $M^{\top}$. $\qs$ 
\enle

By Proposition \ref{p:PropG}, the quotient 
$T_M:=(\mathbb{H}\times\mathbb{C}^n)/G_M$ is 
a complex manifold of complex dimension $n+1$. 
If $n=1$, the manifold $T_M$ is called  {\it Inoue surface} 
(see \cite{Inoue1974}). 
Since the action of $H_M$ on $\mathbb{H}\times\mathbb{C}^n$ 
is also free and properly discontinuous, 
$C_M:=(\mathbb{H}\times\mathbb{C}^n)/H_M$ is also 
a complex manifold of dimension $n+1$. 
If we consider $\mathbb{H}$ as 
$\sqrt{-1}\mathbb{R^*_+}\times\mathbb{R}$, 
the group $H_M$ acts on $\sqrt{-1}\mathbb{R^*_+}$ trivially. 
The quotient $(\mathbb{R}\times\mathbb{C}^n)/H_M$ 
is a $(2n+1)$-dimensional torus $\ttt^{2n+1}$.
The map $g_M$ descends to a self-diffeomorphism of $\mathbb{T}^{2n+1}$. 
Thus we have $C_M=\sqrt{-1}\mathbb{R^*_+}\times\mathbb{T}^{2n+1}$. 
Observe that the matrix $M^\top$ determines a 
self-diffeomorphism of $\ttt^{2n+1}$, this diffeomorphism will be 
denoted by the same symbol 
$M^\top$.

\bigskip

\noindent
\bepr\lb{p:PropI} \ 
The  manifold $T_M$ is diffeomorphic to the mapping torus of 
$$
M^\top
:
\ttt^{2n+1}\to \ttt^{2n+1}.
$$
In particular, $T_M$ is compact. 
\enpr
{\it Proof.} \ 
From \ref{p:PropE}, we have the equality 
\[
T_M=(\mathbb{H}\times\mathbb{C}^n)/G_M
=C_M/\langle g_0\rangle
=(\sqrt{-1}\mathbb{R^*_+}\times\mathbb{T}^{2n+1})/\langle g_0\rangle.
\]
This last manifold is diffeomorphic to 
the manifold obtained from $[1,\alpha]\times\mathbb{T}^{2n+1}$ 
by glueing $\{1\}\times\mathbb{T}^{2n+1}$ with 
$\{\alpha\}\times\mathbb{T}^{2n+1}$ by $g_M$. 
The conclusion now follows from Lemma \ref{l:LemmaH}. 
$\qs$

\bigskip

\section{Topological properties of $T_M$}
\lb{s:topo_pr}


We begin by computation of the first Betti number of $T_M$.
Then  we show that the 
homeomorphism type of $T_M$ determines the 
matrix $M$ up to conjugacy in $\ssll(2n+1,\zz)$
and inverting $M$ (see Theorem \ref{t:a_b}).
This result implies in particular 
(Subsection \ref{su:non-diag})
that if $M$ is not diagonalizable,
then the manifold $T_M$ is not homeomorphic
to any of the manifolds constructed in \cite{OT}.

\subsection{The first Betti number}
\lb{su:first-b}
\noindent
\bele\lb{l:LemmaK}
The first Betti number $b_1(T_M)$ of $T_M$ is equal to $1$. 
\enle
\bigskip

\noindent
{\it Proof.} \ 
The fundamental group $\pi_1(T_M)$ of $T_M$ is isomorphic to $G_M$, 
which has the finite presentation given in Corollary \ref{c:CorollaryF}. 
Hence the first homology group $H_1(T_M;\mathbb{Z})$ is isomorphic to 
the abelian group generated by $g_0,g_1,\ldots ,g_{2n+1}$ with relations 
\[
g_i=m_{i1}g_1+\cdots +m_{i,2n+1}g_{2n+1} 
\quad 
(i\in\{1,\ldots ,2n+1\}).
\]
Since $1$ is not an eigenvalue of $M$, 
we have  $\det(M-I)\ne 0$. 
Thus the first homology group $H_1(T_M;\mathbb{Q})$ with rational coefficient 
is isomorphic to $\mathbb{Q}$. 
$\qs$

\subsection{On fundamental groups of mapping tori}
\lb{su:hom_tori}

Let $k$ be a natural number; any matrix $A\in \ssll(k,\zz)$;
determines  a homeomorphism
$\phi_A:\ttt^k\to\ttt^k$.
Denote by $\TTTT_A$ the mapping torus 
of this map,
we have a fibration $p_A:\TTTT_A\to S^1$
with fiber $\ttt^k$.

\beth\lb{t:a_b}
Let $A,B\in \ssll(k,\zz)$, assume that $1$
is not an eigenvalue of $A$ neither of $B$.
Assume that $\pi_1(\TTTT_A)\approx \pi_1(\TTTT_B)$.
Then $A$ is conjugate to $B$ or to $B^{-1}$ in
$\ssll(k,\zz)$.
\end{theo}
\Prf
Consider the infinite cyclic covering 
$\ove{\TTTT_A} \to \TTTT_A$
induced from the universal  covering $\rr\to S^1$
by $p_A$.
The space $\ove{\TTTT_A}$
is homotopy equivalent to 
the fiber    of $p_A$, that is, to $\ttt^k$.
Therefore the 
Milnor exact sequence 
\cite{Milnor}
of the covering $\ove{\TTTT_A} \to \TTTT_A$
is isomorphic to the following sequence
$$
H_1(\ttt^k)\arrr {A-1} 
H_1(\ttt^k)
\arr
H_1(\TTTT_A)
\arr
H_0(\ttt^k)
\arrr {0}
H_0(\ttt^k).
$$
Since $A-1$ is injective,
the group
$H_1(\TTTT_A)$
is isomorphic to $\zz\oplus F$
where $F$ is a finite abelian group.
Therefore there are exactly
two epimorphisms
$\pi_1(\TTTT_A)$ onto $\zz$,
and they are obtained from one another
via multiplication by $(-1)$.
Consider the exact sequence of the fibration
$p_A$:
$$
0\arr 
\pi_1(\ttt^k)
\arr
\pi_1(\TTTT_A)
\arrr {(p_A)_*}
\pi_1(S^1)
\arr 
0.
$$
It follows from this sequence
that 
$\pi_1(\TTTT_A)$
is isomorphic 
to the semi-direct
product
$\zz\ltimes\zz^k$
where the action of the generator
$t$ of $\pi_1(S^1)$
equals $A$.
Denote by $\iota$ the canonical
generator of 
$\pi_1(S^1)$
and choose an element
$\t_A
\in
\pi_1(\TTTT_A)$
\sut~ 
$(p_A)_*(\t_A)=\iota$.
Let $f:
\pi_1(\TTTT_A)
\to
\pi_1(\TTTT_B)
$
be an isomorphism.
It follows from
the remark above 
that the following diagram
is commutative
$$
\xymatrix{ 
\pi_1(\TTTT_A) \ar[d]^f  \ar[r]^{(p_A)_*} & \pi_1(S^1)  
\ar[d]^{\varepsilon}  \ar[r]  & 0\\
\pi_1(\TTTT_B)  \ar[r]^{(p_B)_*} & \pi_1(S^1)  \ar[r] & 0
}
$$
where $\ve$
equals $1$ or $-1$.
Therefore the element $f(\t_A)$ equals $\t_B\cdot g$
or $(\t_B)^{-1}\cdot g$ with some $g\in \pi_1(\ttt^k)$.
Thus the \ho~
$A$ is conjugate to $B$ or to $B^{-1}$ in 
$\ssll(k, \zz)$.
$\qs$
\beco\lb{c:diago}
If $\pi_1(\TTTT_A)\approx \pi_1(\TTTT_B)$
and $A$ is 
diagonalizable,
then $B$ is also
diagonalizable. $\qs$
\end{coro}
Theorem \ref{t:a_b} above 
can be reformulated in terms of
semi-direct products of groups.
Let $A\in \ssll(k,\zz)$.
Consider the action $\circ$ of $\zz$ on $\zz^k$ defined by 
$m\circ x=A^m\cdot x$;
denote by $S_A$ the corresponding 
semi-direct product 
$\zz\ltimes \zz^k$.

\beco\lb{c:reform-semidir}
Let $A,B\in \ssll(k,\zz)$, assume that $1$
is not an eigenvalue of $A$ neither of $B$.
Assume that $S_A\approx S_B$.
Then $A$ is conjugate to $B$ or to $B^{-1}$ in
$\ssll(k,\zz)$.
\end{coro}
\Prf 
Define an action $\circ$ of $S_A$ on $\rr\times\rr^k$ as follows:
$$
(m,h)\circ(t,v)
=
(t+m, A^mv+h)
$$
The quotient space is clearly $K(S_A,1)$.
It is also easy to see that it is the mapping torus
of the map $A:\ttt^k\to\ttt^k$.
Thus $S_A\approx S_B$
implies 
$\pi_1(\TTTT_A)\approx \pi_1(\TTTT_B)$,
and applying the preceding theorem we deduce the Corollary. $\qs$

\bede\lb{d:diag-type}
We say that the semi-direct product $S_A=\zz\ltimes\zz^k$
is {\it of diagonal type}, if 
$A$ is diagonalizable over $\cc$
and its eigenvalues are different from 1.


We say that the semi-direct product $S_A=\zz\ltimes\zz^k$
is {\it of non-diagonal type}, if
$A$ is non-diagonalizable over $\cc$
and its eigenvalues are different from 1.


\end{defi}

\beco\lb{c:diag-type}
A semi-direct product of diagonal type
is not isomorphic to a semi-direct product 
of non-diagonal type. $\qs$
\enco

\bepr\lb{p:orient}
Let $A,B\in \ssll(2n+1,\zz)$.
The \ma s $\TTTT_A$ and $\TTTT_B$
have then natural   orientations.
Assume that $A$ is conjugate to $B$ or to $B^{-1} $
in $\ssll(2n+1,\zz)$.
Then there is an orientation preserving 
diffeomorphism $\TTTT_A\approx \TTTT_B$.
\enpr
\Prf
1) 
If $A= C^{-1}BC$ with 
$C\in \ssll(\zz,2n+1)$
then the required diffeomorphism
is given  by the formula
$(x,t)\mapsto (Cx, t)$.

2) 
If $A= C^{-1}B^{-1}C$ with 
$C\in \ssll(\zz,2n+1)$
then the required diffeomorphism
is defined
as the composition
$\chi\circ \phi\circ\psi$
with 
\begin{gather*}
 \psi : \TTTT_A \to \TTTT_A, \ \ \psi(x,t) = (-x,t), \\
  \chi : \TTTT_{B^{-1}} \to \TTTT_B, \ \ \chi(x,t) = (x,1-t), \\
  \phi: \TTTT_A \to \TTTT_{B^{-1}},  \ \ \phi(x,t) = (Cx,t) 
\end{gather*}
Observe that $\psi$ and $\chi$ reverse  orientation,
and $\phi$ is orientation preserving. 
The proposition is proved. $\qs$

\section{Geometric properties of $T_M$}
\label{s:geom_pr}

This section is about the properties
of the manifolds $T_M$ related to its complex structure.
These properties are mostly similar to the properties of OT-manifolds.
The first section is about the holomorphic bundles over $T_M$ and their sections.

In the last two subsections we investigate the questions of existence 
of  K\"ahler and locally conformally K\"ahler structures on manifolds $T_M$.
Here we concentrate ourselves on the case
when the matrix $M$ is not diagonalizable.

\subsection{Holomorphic bundles on $T_M$ and their sections}
\lb{su:holom}

\bigskip 

\noindent
\bepr\lb{p:holomm} 
Any $H_M$-invariant holomorphic function on $\mathbb{H}\times\mathbb{C}^n$ 
is constant. 
\enpr
\bigskip

\noindent
{\it Proof.} \ 
Let $f:\mathbb{H}\times\mathbb{C}^n\rightarrow \mathbb{C}$ be a holomorphic function. 
Suppose that $f(g(w,z))=f(w,z)$ for every $g\in H_M$ and 
$(w,z)\in \mathbb{H}\times\mathbb{C}^n$. 
In particular, for $i\geq 1$ we have 
\bq\lb{f:shift}
f(w,z)=f(g_i(w,z))=
f((w,z)+u_i).
\end{equation}
For $w_0\in\hh$ let
$$
A_w=\{(w_0,z)~|~ z\in\cc^n\}, \ \ B_w=\{(w,z)~|~ z\in\cc^n, \ {\frak Im\ }  w = {\frak Im\ } w_o\}.
$$
Then  $A_w$ is an $n$-dimensional complex space, and $B_w$
is a real vector space of dimension $2n+1$. We have $A_w\sbs B_w$.
The abelian group generated by the vectors $u_1, \ldots, v_{2n+1}$
is a full lattice in $B_w$, therefore $f~|~B_w$ is bounded, and so is 
$f~|~A_w$. The function $f~|~A_w$ is holomorphic and bounded, therefore 
it is constant. Thus $f(w,0)=f(w,z)$ for every $(w,z)\in\HC$.
Consider a subset 
\[
A:=
\left\{ \left. 
\sum_{i=1}^{2n+1} s_ia^{(i)} \, \right|\, 
s_1,\ldots ,s_{2n+1}\in\mathbb{Z}
\right\} \sbs \rr.
\]
Using \rrf{f:shift} repeatedly, we deduce that
$f(w,0)=f(w+\xi,0)$ 
for every $\xi\in A$.
Since $a$ is an eigenvector of $M$ 
corresponding to $\alpha$, 
we have $\alpha a^{(i)}\in A$ for every $i$. 
The set 
$A_0:=\{(n_1+n_2\alpha)a^{(i)} \, |\, n_1,n_2\in\mathbb{Z}\}$ is included in $A$, 
and it is dense in $\mathbb{R}$ by Kronecker's density theorem. 
Therefore $A$ is also dense in $\mathbb{R}$,
and $f(w,0)$ does not depend on $w$. $\qs$ 

\pa
Our next proposition is similar to \cite{OT}, Prop. 2.5.
\bepr\lb{p:holo_sects}
\been\item
There are no non-trivial holomorphic 1-forms on $T_M$.
\item
Let $K=K_{T_M}$ be the canonical bundle on $T_M$.
Let $k\in \nn, k>0$. Then the bundle $K^{\otimes k}$
admits no non-trivial global  sections. 
The Kodaira dimension of $T_M$ is therefore equal to $-\infty$.
\enen
\enpr

\Prf
1) Let $\l$ be a holomorphic 1-form on $T_M$.
Denote by $u:\hh\times \cc^n\to T_M$ the universal covering of $T_M$.
We have
$$
u^*\l=f_0(w,z)dw+f_1(w,z)dz_1+ \ldots + f_1(w,z)dz_n,
$$
where $f_i$ are holomorphic functions on $\hh\times \cc^n$;
they are invariant \wrt~ $H_M$, therefore constant by 
Proposition \ref{p:holomm}.
The form $u^*\l$ is also $g_0$-invariant. Since $g_0^*(dw)=\a dw$, 
we have $f_0=0$. Similarly, since 1 is not an eigenvalue of $M$ we deduce
that $f_i(w,z)=0$ for every $i\geq 1$. 

2) Let $\r$ be a section of 
$K^{\otimes k}$.
Then
$$u^*\r
=
f(w,z)\Big(dw\wedge dz_1 \wedge\ldots \wedge dz_n\Big).
$$
Similarly to the point 1) we deduce that $f(z,w)$ is a constant function.
Since $u^*\r$ is also $g_0$-invariant, we have
$(\a \cdot \b_1\cdots \b_n)^k = 1$. The condition $\det M=1$ implies then that 
 $(\bar\b_1\cdots \bar\b_n)^k = 1$, and finally $\a^k=1$, which is impossible. $\qs$

\subsection{K\"ahler structures}
\label{su:formality}

It is clear that the manifold $T_M$ 
is not 
K\"ahler,
since $b_1(T_M)=1$.
A similar argument implies that  the manifold $(T_M)^l$ is not 
 K\"ahler if $l$ is odd.
  The next proposition extends this property to all $l\in\nn$.
 \bepr\lb{p:nonkel}
 For any $l\in\nn$ the complex manifold $X=(T_M)^l$
does not admit a 
 K\"ahler
structure.
\enpr
\Prf
Consider the composition 
$\pi':
(T_M)^l
\arrr {p_1}
T_M
\arrr {\pi}
S^1
$
where $\pi$ is the fibration 
induced by the mapping torus 
structure on $T_M$.
The map $\pi'$ is a fibration
with fiber $(T_M)^{l-1}\times \ttt^{2n+1}$.
The monodromy \ho~ of this fibration 
equals $\Id\times M^\top$.
This matrix is not diagonalizable,
and the main theorem of \cite{Pjord1}
implies that $(T_M)^l$
does not admit a 
 K\"ahler
structure.
$\qs$ 

\subsection{Locally conformally K\"ahler
structures}
\lb{su:lck}

In 1982 F. Tricerri \cite{Tr}
proved that Inoue manifold admits an LCK-structure.
The case of OT-manifolds is different,
it is proved in \cite{OT}
that the OT-manifolds $X(K,U)$
do not admit an LCK-structure for $s=1$.

\bepr\lb{p:nonLCK}
Assume that $M$ is not 
diagonalizable.
Then $T_M$ does not admit 
an LCK-structure.
\enpr
\Prf 
The proof follows the lines of 
the corresponding theorem of Oeljeklaus-Toma
\cite{OT}, Prop. 2.9; in our case the argument 
is somewhat simpler.
Assume that there exists an LCK-structure 
on $T_M$ ;
denote by $\Omega$ the corresponding 1-form,
so that $d\Omega=\o\wedge \O$. 
Consider the infinite cyclic covering
$p:\ove{T_M}\to T_M$
corresponding to the mapping torus structure 
of $T_M$. The universal covering 
$
\hh\times \cc^n
\to T_M
$ factors
as follows:
$$
\hh\times \cc^n
\arrr q
\ove{T_M}
\arrr p T_M.
$$
We have a diffeomorphism
$
\ove{T_M}\approx \ttt^{2n+1} \times \rr$,
and replacing 
the form $(q\circ p)^*\O$
by its average \wrt~ the action of $\ttt^{2n+1}$ 
we can assume that the form 
$(q\circ p)^*\O$ on $\hh\times \cc^n$
does not depend on the coordinates $z=(z_1, \ldots , z_n)$
on every subspace $\{h\}\times \cc^n$.
Let $(q\circ p)^*\O= df$, with $f:\hh\times \cc^n\to \rr$.
Since $(q\circ p)^*\O$ is a symplectic form on 
$\{h\}\times \cc^n$,
this implies $df=0$
on $\{h\}\times \cc^n$.
Put 
$\tau= e^{-f}\cdot (q\circ p)^*\O$. Then
$$
\tau
=
\sum_{0\leq i < j \leq n} g_{ij}(z)dz_i\wedge d \bar z_j,
$$
(where we have named the first coordinate
of 
$\{h\}\times \cc^n$
by $z_0$)
Here $g_{ij}(z_0, z_1, \cdot, z_n)$
does not depend on 
$(z_1, \cdot, z_n)$.
Moreover, $d\tau=0$, and this implies easily that
$g_{ij}(z)$ does not depend on $z_0$ either.
We can assume that $f(\sqrt{-1}, 0, \ldots, 0)=0$.
Let $\xi=
f(\sqrt{-1}, 0, \ldots, 0)=0\in\rr$, and $\mu=e^{-\xi}$.
Denote by $\tau_0$ the restriction of $\tau$ to $i\times \cc^n$.
Then we have
$(M^\top)^*\tau_0=\mu\cdot\tau_0$, which implies
that the linear map $M^\top/\mu$ 
conserves the non-degenerate 2-form $\tau_0\in\wedge^2(\cc^n)$.
The symmetric form
$\s(x,y)=\tau_0(x, iy)$ on $\cc^n$
is a scalar product (since $\O$ is the imaginary part of a hermitian form), therefore $M^\top/\mu$
conserves a  scalar product, which is impossible since 
$M^\top$ is non-diagonalizable. $\qs$

\section{Relations with the Oeljeklaus-Toma construction}
\lb{s:Tm-OT}


In this section we study the relation 
between the \ma~ $T_M$ constructed in \S \ref{s:constr}
and the \ma s constructed by \OT~ in \cite{OT}
({\it OT-manifolds} for short). In Subsection 
\ref{su:ot-as-tm} we show that some of 
OT-manifolds appear as  $T_M$-manifolds.
In Subsection \ref{su:non-diag}
we show that the manifold $T_M$ with $M$ non-diagonalizable 
is not homeomorphic to any of OT-manifolds.

\subsection{Construction of  OT-manifolds }
\lb{su:constrOT}

Let us first recall the construction from \cite{OT}
(in a slightly modified  terminology).
Let $K$ be an algebraic number field.
An embedding $K \hookrightarrow \cc$ 
is called {\it real} if its image is     in $\rr$;
an embedding which is not real is called {\it complex}.
Denote by $s$ the number  of real embeddings 
and by $t$ the number  of complex embeddings.
Then $(K:\qq)=s+2t$.
Let $\s_1, \ldots, \s_s$ be the real embeddings
and $\s_{s+1}, \ldots, \s_{s+2t}$
be the complex  embeddings;
we can assume             that 
$\s_i=\ove{\s_{t+i}}$ for $i\geq s+1$.
The map 
$$\s:K\to\rr^s\times \cc^t;
\ \ \ 
\s(x)=\big(\s_1(x), \ldots , \s_{s+t}(x)\big)
$$
is an embedding (known as {\it geometric representation}
of the field $K$, see \cite{BS}, Ch. II, \S 3).
Let $\OOOO$ be any order in $K$, then 
$\s(\OOOO)$ is a full lattice in 
$
\rr^s\times \cc^r.
$
Denote by $\OOOO^*$ the group of all units of $\OOOO$.
Dirichlet theorem (see \cite{BS}, Ch. II, \S 4 , Th. 5)
says that the group 
$\OOOO^*/Tors$
is a free abelian group of rank $s+t-1$.
Assume that $t\geq 1$.
Choose any elements  $u_1, \ldots, u_s$ 
of this group generating a free abelian subgroup
of rank $s$. A unit $\l\in\OOOO$ 
will be called {\it positive}
if $\s_i(\l)>0$ for every $i\leq s$.
Replacing
$u_i$ by $u_i^2$ if necessary 
we can assume that 
every $u_i$ is positive.
The subgroup $U$ of 
$\OOOO^*/Tors$
generated by $u_1, \ldots, u_s$ 
acts on $\OOOO$
and we can form 
the semi-direct product $\PPPP=U\ltimes \OOOO$.
The group $\PPPP$
acts on 
$\cc^r=\cc^s\times \cc^t$
as follows:
\begin{itemize}
 \item 
 any element $\xi\in\OOOO$
 acts by translation by vector
 $\s(\xi)\in \rr^s\times \cc^r.$
 \item any element $\l\in U$ acts
 as follows:
 $$
\l\cdot (z_1, \ldots, z_{s+t})
=
(\s_1(\l)z_1, \ldots, \s_{s+t} (\l)z_{s+t}).
$$
\end{itemize}
For $i\leq s$ the numbers $\s_i(\l) $
are  real and positive, so the subset
$\hh^s\times\cc^t$ is invariant under the action
of $\PPPP$.
This action is properly discontinuous
and the quotient is a complex analytic
manifold of dimension $s+t$
which will be denoted by $X(K,\OOOO, U)$.
The notation $X(K,U)$ used in the article \cite{OT}
pertains to the case when the order $\OOOO$ 
is the maximal order of $K$.

\subsection{OT-manifolds  as manifolds of type $T_M$}
\lb{su:ot-as-tm}

Consider the case $s=1$. In this subsection
we will denote the number of complex embeddings
of $K$ by $n$, in order
to fit to the terminology of the previous sections.
Then $(K:\qq)=2n+1$. We assume that $n\geq 1$.
Assume that there is a Dirichlet unit $\xi$ in $K$ 
\sut~ $\qq(\xi)=K$.
This assumption 
holds for example 
when there are no proper subfields 
$\qq\subsetneq K'\subsetneq K$; this is always the 
case if $2n+1$ is a prime number.
Replacing $\xi$ by $\xi^2$ if necessary
we can assume that $\xi$ is positive.
Denote by $\OOOO$ the order $\zz[\xi]$, and let $U$ be the group 
of units, generated by $\xi$.
Denote by $P$ the minimal polynomial of $\xi$,
let $C_P$ be the companion matrix of $P$,
and put $D_P=C_P^\top$.

\bepr\lb{p:ot-tm}
We have a biholomorphism
$$T_{D_P}\approx X(K,\OOOO, U).$$
\enpr
\Prf
Let us give explicit descriptions of both these \ma s.
\pa
{\bf 1. }  \  The \ma~ $X(K,\OOOO, U).$
\pa
\noindent
The lattice $\s(\OOOO)$ is a free $\zz$-module generated
by $e_i=\s(\xi^i)$.
Denote by $\a, \b_1, \ldots, \b_n, \bar\b_1, \ldots, \bar\b_n,$
the roots of $P$ (here $\a\in\rr, \b_i\notin \rr$).
Then 
$\s(\xi^k)= (\a^k, \b^k_1, \ldots, \b^k_n)$, and the action of 
$\xi$ on this basis is given by 
the following formula:
$$
\xi\cdot (x,z_1, \ldots, z_n)
=
(\a x, \b_1z_1  , \ldots,  \b_n z_n).
$$
\pa
{\bf 2. }  \  The \ma~ $T_{D_P}.$
\pa
\noindent
The eigenvalues of the matrix $D_P$
are the same as of the matrix $C_P$, that is,
$\a, \b_1, \ldots, \b_n, \bar\b_1, \ldots, \bar\b_n,$.
The corresponding eigenvectors are:
$$
a=(1,\a,\ldots, \a^{2n}), \ \ 
b_i=(1,\b_i,\ldots, \b_i^{2n})
\
({\rm where \ }
1 \leq i \leq n).
$$
The vectors $u_i$
generating the group $H_{D_P}$ of translations 
(see Section \ref{s:constr}, page \pageref{f:M})
are 
given by the formula
$$
u_1=(1, \ldots, 1),
\
u_2=(\a, \b_1, \ldots, \b_n), 
\ 
\ldots 
\ 
,
u_{2n+1}=(\a^{2n}, \b_1^{2n}, \ldots, \b_n^{2n}).
$$
Both  matrices $D_{P}, \ C_P$
restricted to the subspace $W$ generated by 
the vectors $u_i$ with $2\leq i \leq n+1$
are diagonal in this bases
and the diagonal entries are equal to 
$\b_1, \ldots, \b_n$.
Therefore the element $g_0\in G_{D_P}$ 
acts as follows
$$
g_0\cdot (w,z_1, \ldots, z_n)
=
(\a w, \b_1z_1  , \ldots,  \b_n z_n).
$$

The proposition follows. $\qs$

\subsection{The case of non-diagonalizable matrix $M$}
\lb{su:non-diag}

\bele\lb{l:OT-diag-type}
Let $K$ be an algebraic number field with $s=1$, put $(K:\qq)=2n+1$.
Let $\OOOO$ be an order in $K$, and $\xi$ a positive unit of $\OOOO$.
Let $X=X(K,\OOOO, \xi)$ be the corresponding OT-manifold.
Then the group $\pi_1(X)$ is a semi-direct product 
$\zz\ltimes\zz^{2n+1}$ of diagonal type.
\enle
\Prf
The $\pi_1(X)$ is a semi-direct product $S_A$ 
where $A$ is the matrix of the action  of the unit $\xi$
on $\OOOO$. Let $P$ be the minimal polynomial of $\xi$.
The roots of $P$ are simple and different from 1. 
Since $P(A)=0$, the minimal polynomial of $A$ has the same properties.
Therefore $A$ is diagonal and Lemma is proved. $\qs$

\bepr\lb{p:compar}
Let $M\in \ssll(2n+1,\zz)$
be a matrix of type $\II$ and non-diagonalizable over $\cc$.
Then the group $\pi_1(T_M)$
is not isomorphic to 
the fundamental group of
any of manifolds
$X(K,U)$ constructed in \cite{OT}.
Therefore $T_M$ is not homeomorphic to any of manifolds
$X(K,U)$.
\enpr
\Prf
Assume that we have an isomorphism
$\pi_1(X(K,U))\approx \pi_1(T_M)$.
Since $b_1(\pi_1(X(K,U)))=s$, and $b_1(\pi_1(T_M))=1$, we have $s=1$. 
By Lemma \ref{l:OT-diag-type} $\pi_1(X(K,U))$
is isomorphic to a semi-direct product $\zz\ltimes\zz^{2n+1}$
of diagonal type. Recall that  $\pi_1(T_M)$
is a a semi-direct product
of a {\it non-diagonal type}.
Apply Corollary  \ref{c:diag-type}
and the proof is over.
$\qs$ 

\section{Acknowledgements}
\lb{s:ack}

The second author thanks Professor F. Bogomolov for 
initiating him to the theory of Inoue surfaces in 2013,
for several discussions on this subject and for support.
The work on this article began in January 2018 
when the second author was visiting the Tokyo Institute of
Technology; many thanks for the warm hospitality and support!

\end{document}